\documentclass[12pt]{article}

\usepackage[margin=1in]{geometry}  
\usepackage{graphicx}              
\usepackage{amsmath}               
\usepackage{amsfonts}              
\usepackage{amsthm}                
\usepackage{amssymb}

\begin{document}

\nocite{*}

\title{Hermitian forms of K3 type.}

\author{Evgeny Mayanskiy}

\maketitle

\begin{abstract}
We classify quadratic spaces over endomorphism fields of $K3$ surfaces. We consider both totally real and CM cases. 
\end{abstract}



\medskip

Let $E$ be a number field, $[E\colon \mathbb Q]=r+2s$, where $r=\mid {\Sigma}_{\mathbb R}\mid $ is the number of real embeddings of $E$ (${\Sigma}_{\mathbb R} = \{  \sigma \colon E\rightarrow \mathbb R \}$) and $s=\mid {\Sigma}_{\mathbb C}\mid $ is the number of pairs (of complex-conjugate) complex embeddings of $E$ (${\Sigma}_{\mathbb C} = \{  ( \sigma \colon E\rightarrow \mathbb C, \bar{\sigma} \colon E\rightarrow \mathbb C ) \; \mid \; {\sigma}\neq {\bar{\sigma}} \}$).\\

Let us assume that either $E$ is totally real (and so $s=0$) or $E=E_0(\theta)$ is a CM-field (and so $r=0$), where $E_0$ is a totally real field (fixed by the CM-involution) and ${\theta}^2\in E_0$, $\theta \notin E_0$.\\

Let $\{ e_1, \cdots , e_{r+2s}  \}$ be a basis of $E / \mathbb Q$ (and assume that $e_{r+s+i}=\theta \cdot e_{r+i}$ for any $i\geq 1$). Denote by $\xi \mapsto \bar{\xi}$ the CM-involution on $E$ (which is the identity, if $E$ is totally real).\\

Let $V$ be a vector space over $E$ of dimension $m=dim_{E}V\geq 1$. Let $\Phi \colon V {\otimes}_E V\rightarrow E$ be a nondegenerate hermitian sesquilinear (with respect to the CM-involution on $E$) form on $V$.\\

Let $\{ u_1, \cdots , u_m \}$ be a basis of $V$ over $E$.\\

{\bf Lemma 1.} {\it If $E$ is totally real, then there is an isomorphism of $\mathbb R$-vector spaces:
$$
\gamma \colon V{\otimes}_{\mathbb Q} \mathbb R \xrightarrow{~} {\bigoplus}_{\sigma \in {\Sigma}_{\mathbb R}} V_{\sigma},\quad \quad (v,c)\mapsto ({\sigma}(v)\cdot c)_{\sigma}.
$$
This decomposition is orthogonal with respect to the symmetric bilinear form $tr\circ \Phi \colon V{\otimes}_{\mathbb Q} V\rightarrow \mathbb Q$ (where $tr \colon E\rightarrow \mathbb Q$ denotes the trace morphism).\\

For any $\sigma \in {\Sigma}_{\mathbb R}$ the restriction of $tr\circ \Phi$ onto $V_{\sigma}$ is the localization of $\Phi$ at $\sigma$: 
$$
{\Phi}_{\sigma}\colon {V_{\sigma}}{\otimes}_{\mathbb R} {V_{\sigma}}\rightarrow \mathbb R, \quad \quad (\xi, \zeta)\mapsto {\sigma}({\Phi}(\xi,\zeta))
$$}

{\it Proof:} Let $V_{\sigma}={\gamma}^{-1}(V_{\sigma})=$ $\{ \xi \in V {\otimes}_{\mathbb Q} \mathbb R \; \mid \; e\cdot \xi = \xi \cdot {\sigma}(e) $ for any $ e\in E \}=$ $V {\otimes}_{E, \sigma} \mathbb R$.\\

For any $\sigma \in {\Sigma}_{\mathbb R}$ and $i=1,...,m$ choose constants $c_{i}^{(\sigma)}\in \mathbb R$ such that $\sum_{i}{c_{i}^{(\sigma)}\cdot {{\sigma}'}(e_i)}={\delta}_{\sigma,{\sigma}'}$ (Kroneker symbol) for any $\sigma, {\sigma}'\in {\Sigma}_{\mathbb R}$. Then for any $\sigma$ and $i$ we have: ${\gamma}^{-1}(u_i {\otimes}_{E,\sigma} 1)=\sum_{j}{(e_j\cdot u_i){\otimes}_{\mathbb Q} {c_{j}^{(\sigma)}}}$.\\ 

For any $\sigma, {\sigma}'\in {\Sigma}_{\mathbb R}$, $(tr\circ \Phi)(u_i {\otimes}_{E,\sigma} 1, u_j {\otimes}_{E,{\sigma}'} 1)=$ $\sum_{k,l}{(tr\circ \Phi)( (e_k\cdot u_i){\otimes}_{\mathbb Q} {c_{k}^{(\sigma)}} , (e_l\cdot u_j){\otimes}_{\mathbb Q} {c_{l}^{({\sigma}')}} )}=$ $\sum_{k,l}{c_{k}^{(\sigma)}\cdot c_{l}^{({\sigma}')} \cdot tr (e_ke_l\cdot {\Phi}(u_i,u_j)) }=$ $\sum_{{\sigma}''}{\sum_{k,l}{  c_{k}^{(\sigma)} c_{l}^{({\sigma}')} \cdot {{{\sigma}''}(e_k)}\cdot {{{\sigma}''}(e_l)} \cdot {{{\sigma}''}({\Phi}(u_i,u_j))}  }}=$ ${\delta}_{\sigma,{\sigma}'}\cdot {\sigma}({\Phi}(u_i,u_j))$. {\it QED}\\

{\bf Lemma 1'.} {\it If $E$ is a CM-field, then there is an isomorphism of $\mathbb C$-vector spaces:
$$
\gamma \colon V{\otimes}_{\mathbb Q} \mathbb C \xrightarrow{~} {\bigoplus}_{\sigma \in {\Sigma}_{\mathbb C}} (V_{\sigma} \oplus V_{\bar{\sigma}} ),\quad \quad (v,c)\mapsto ({\sigma}(v)\cdot c, {\bar{\sigma}}(v)\cdot c)_{\sigma}.
$$
This decomposition (over $\sigma \in {\Sigma}_{\mathbb C}$) is orthogonal with respect to the symmetric bilinear form $tr\circ \Phi \colon V{\otimes}_{\mathbb Q} V\rightarrow \mathbb Q$.\\

For any $\sigma \in {\Sigma}_{\mathbb C}$ the restriction of $tr\circ \Phi$ onto $V_{\sigma}$ is the localization of $\Phi$ at $\sigma$: $(tr\circ \Phi)(V_{\sigma},V_{\sigma})=0$ and 
$$
(tr\circ \Phi){\mid}_{V_{\sigma}}={\Phi}_{\sigma}\colon {V_{\sigma}}{\otimes}_{\mathbb C} {V_{\bar{\sigma}}}\rightarrow \mathbb C, \quad \quad (\xi, \bar{\zeta})\mapsto {\sigma}({\Phi}(\xi,\zeta))
$$}

{\it Proof:} Let $V_{\sigma}={\gamma}^{-1}(V_{\sigma})=$ $\{ \xi \in V {\otimes}_{\mathbb Q} \mathbb C \; \mid \; e\cdot \xi = \xi \cdot {\sigma}(e) $ for any $ e\in E \}=$ $V {\otimes}_{E, \sigma} \mathbb C$.\\

For any $\sigma \in {\Sigma}_{\mathbb C}$ and $i=1,...,m$ choose constants $c_{i}^{(\sigma)}\in \mathbb C$ such that $\sum_{i}{c_{i}^{(\sigma)}\cdot {{\sigma}'}(e_i)}={\delta}_{\sigma,{\sigma}'}$ for any $\sigma, {\sigma}'\in {\Sigma}_{\mathbb C}$. Then for any $\sigma$ and $i$ we have: ${\gamma}^{-1}(u_i {\otimes}_{E,\sigma} 1)=\sum_{j}{(e_j\cdot u_i){\otimes}_{\mathbb Q} {c_{j}^{(\sigma)}}}$.\\ 

Note that $\overline{c_{i}^{(\sigma)}}=c_{i}^{(\bar{\sigma})}$ and $\overline{V_{\sigma}}=V_{\bar{\sigma}}$ as $\mathbb C$-vector subspaces of $V {\otimes}_{\mathbb Q}\mathbb C$.\\

For any $\sigma, {\sigma}'\in {\Sigma}_{\mathbb C}$, $(tr\circ \Phi)(u_i {\otimes}_{E,\sigma} 1, u_j {\otimes}_{E,{\sigma}'} 1)=$ ${\delta}_{\sigma,\overline{{\sigma}'}}\cdot {\sigma}({\Phi}(u_i,u_j))$ by the same computation as in the previous Lemma. {\it QED}\\

{\bf Corollary 1.} {\it $tr\circ \Phi \colon V{\otimes}_{\mathbb Q} V\rightarrow \mathbb Q$ is a nondegenerate symmetric bilinear form.}\\

{\bf Corollary 2.} {\it If $E$ is a CM-field, then there is an orthogonal (with respect to $tr\circ \Phi$) decomposition:
$$
V {\otimes}_{\mathbb Q} \mathbb R = {\bigoplus}_{\sigma \in {\Sigma}_{\mathbb C}} (W_{\sigma}\oplus W_{\sigma})
$$
such that the restriction of $tr\circ \Phi$ onto $W_{\sigma}$ is the localization of the nondegenerate symmetric $E_0$-bilinear form ${\Phi}_0\colon V {\otimes}_{E_0} V\rightarrow E_0$, $(\xi,\zeta)\mapsto {\Phi}(\xi,\zeta)+{\Phi}(\zeta,\xi)$ induced by $\Phi$ at the embedding ${\sigma}{\mid}_{E_0}\colon E_0\rightarrow \mathbb R$ of $E_0$ corresponding to $\sigma \in {\Sigma}_{\mathbb C}$.}\\

{\it Proof:} Since $\{ u_1 {\otimes}_{E,\sigma} 1, \cdots ,  u_m {\otimes}_{E,\sigma} 1, u_1 {\otimes}_{E,\bar{\sigma}} 1,\cdots u_m {\otimes}_{E,\bar{\sigma}} 1 \}$ is a $\mathbb C$-basis of $V_{\sigma}\oplus V_{\bar{\sigma}}\subset V {\otimes}_{\mathbb Q} \mathbb C$, $\{  u_1 {\otimes}_{E,{\sigma}} 1 + u_1 {\otimes}_{E,\bar{\sigma}} 1, u_1 {\otimes}_{E,{\sigma}} \sqrt{-1} - u_1 {\otimes}_{E,\bar{\sigma}} \sqrt{-1},\cdots , u_m {\otimes}_{E,{\sigma}} 1 + u_m {\otimes}_{E,\bar{\sigma}} 1, u_m {\otimes}_{E,{\sigma}} \sqrt{-1} - u_m {\otimes}_{E,\bar{\sigma}} \sqrt{-1} \}$ is an $\mathbb R$-basis of $(V {\otimes}_{\mathbb Q} \mathbb R) \cap (V_{\sigma}\oplus V_{\bar{\sigma}}) = W_{\sigma}\oplus {W'}_{\sigma}$.\\

We define $W_{\sigma}={\mathbb R} (u_1 {\otimes}_{E,{\sigma}} 1 + u_1 {\otimes}_{E,\bar{\sigma}} 1) + \cdots + {\mathbb R} (u_m {\otimes}_{E,{\sigma}} 1 + u_m {\otimes}_{E,\bar{\sigma}} 1)$ and ${W'}_{\sigma}={\mathbb R} (u_1 {\otimes}_{E,{\sigma}} \sqrt{-1} - u_1 {\otimes}_{E,\bar{\sigma}} \sqrt{-1}) + \cdots + {\mathbb R} (u_m {\otimes}_{E,{\sigma}} \sqrt{-1} - u_m {\otimes}_{E,\bar{\sigma}} \sqrt{-1})$.\\

Note that $(tr\circ \Phi)(u_i {\otimes}_{E,{\sigma}} 1 + u_i {\otimes}_{E,\bar{\sigma}} 1, u_j {\otimes}_{E,{\sigma}} \sqrt{-1} - u_j {\otimes}_{E,\bar{\sigma}} \sqrt{-1})=0$, i.e. $W_{\sigma}$ and ${W'}_{\sigma}$ are orthogonal with respect to $tr\circ \Phi$, $(tr\circ \Phi)(u_i {\otimes}_{E,{\sigma}} 1 + u_i {\otimes}_{E,\bar{\sigma}} 1, u_j {\otimes}_{E,{\sigma}} 1 + u_j {\otimes}_{E,\bar{\sigma}} 1)=$  $ (tr\circ \Phi)(u_i {\otimes}_{E,{\sigma}} \sqrt{-1} - u_i {\otimes}_{E,\bar{\sigma}} \sqrt{-1}, u_j {\otimes}_{E,{\sigma}} \sqrt{-1} - u_j {\otimes}_{E,\bar{\sigma}} \sqrt{-1}) = $  ${\sigma}({\Phi}(u_i, u_j))+{\bar{\sigma}}({\Phi}(u_i,u_j))$. {\it QED}\\

We will say that the $E$-quadratic space $(V,\Phi)$ is {\it of $K3$ type}, if the $\mathbb Q$-vector space $V$ admits a weight $2$ Hodge structure, which is polarized by the bilinear form $tr\circ \Phi$ and coincides with the polarized Hodge structure on the transcendental $\mathbb Q$-lattice of a $K3$ surface and such that there is an isomorphism of $\mathbb Q$-algebras $E\cong End_{Hdg}(V)$ (where $End_{Hdg}(V)$ denotes the algebra of endomorphisms of $V$ preserving the Hodge decomposition), which identifies the structures of a $E$-module and of a $End_{Hdg}(V)$-module on $V$.\\ 

We will say that a nondegenerate hermitian sesquilinear form $\Phi \colon H {\otimes}_{O_E} H\rightarrow O_E$ (where $O_E$ is the ring of integers of $E$) on a free and finitely generated $O_E$-module $H$ is {\it of $K3$ type}, if the abelian group $H$ admits a weight $2$ Hodge structure, which is polarized by the bilinear form $tr\circ \Phi$ and coincides with the polarized Hodge structure on the transcendental lattice of a $K3$ surface and such that there is an isomorphism of $\mathbb Q$-algebras $E\cong End_{Hdg}(H {\otimes}_{\mathbb Z}\mathbb Q)$, which identifies the structures of a $E$-module and of a $End_{Hdg}(H {\otimes}_{\mathbb Z}\mathbb Q)$-module on $H {\otimes}_{\mathbb Z}\mathbb Q$.\\ 

Let $\Lambda = E_8(-1)^{\oplus 2}\oplus U^{\oplus 3}$ be the second cohomology lattice of a $K3$ surface.\\

It is clear that if $(V,\Phi)$ (respectively, $(H,\Phi)$) is of $K3$ type, then there is an embedding of $\mathbb Q$-lattices $(V,tr\circ \Phi)\hookrightarrow \Lambda {\otimes}_{\mathbb Z} \mathbb Q$ (respectively, a primitive embedding of lattices $(H,tr\circ \Phi)\hookrightarrow \Lambda$), as well as the condition on signatures described in the Theorem below. In fact, using the results of Zarhin \cite{Zarhin} one can show that these conditions are also sufficient:\\ 

{\bf Theorem 1.} {\it $(V,\Phi)$ is of $K3$ type, if and only if the following conditions hold:
\begin{enumerate}
\item[(1)] there is an embedding of $\mathbb Q$-lattices $(V,tr\circ \Phi)\hookrightarrow \Lambda {\otimes}_{\mathbb Z} \mathbb Q$,
\item[(2)] if $E$ is totally real, then $m\geq 3$ and there is ${\sigma}_0\in {\Sigma}_{\mathbb R}$ such that $sign ({\Phi}_{{\sigma}_0})=(2,m-2)$ and for any ${\sigma}\in {\Sigma}_{\mathbb R}$, $\sigma \neq {\sigma}_0$, $sign ({\Phi}_{{\sigma}})=(0,m)$, 
\item[(3)] if $E$ is a CM-field, then $m\geq 1$ and there is ${\sigma}_0\in {\Sigma}_{\mathbb C}$ such that $sign ({\Phi}_{{\sigma}_0})=(1,m-1)$ and for any ${\sigma}\in {\Sigma}_{\mathbb C}$, $\sigma \neq {\sigma}_0, \bar{{\sigma}_0}$, $sign ({\Phi}_{{\sigma}})=(0,m)$.
\end{enumerate}}

{\it Proof:} Assuming that these conditions are satisfied, one has to find $f\in V{\otimes}_{\mathbb Q} \mathbb C$ such that $(tr\circ \Phi)(f,f)=0$, $(tr\circ \Phi)(f,\bar{f})>0$ and $ker (V\rightarrow \mathbb C, u\mapsto (tr\circ \Phi)(u,f))=0$. This will give a weight $2$ Hodge structure on $V$ polarized by $tr\circ \Phi$ with $H^{2,0}=\mathbb C f$. Then one needs to make sure that $E=End_{Hdg}(V)$, which in particular means that for any $e\in E$ multiplication by $e$ on $V$ preserves Hodge  decomposition. This implies that $f\in V_{\sigma}\subset V {\otimes}_{\mathbb Q} \mathbb C$ for some $\sigma \in \Sigma = {\Sigma}_{\mathbb R}\cup {\Sigma}_{\mathbb C}$. Moreover, it is clear that one should have $\sigma={\sigma}_0$.\\ 

Let us check that a (very) general $f\in V_{{\sigma}_0}$ works. Let us assume  that $\{ u_1, \cdots , u_m \}$ is an orthogonal basis of $V$, $d_i={\Phi}(u_i,u_i)$ and ${{\sigma}_0}(d_1)\geq {{\sigma}_0}(d_2) \geq \cdots \geq {{\sigma}_0}(d_m)$.\\

Case 1: $E$ is totally real and $m\geq 4$. Let $f=f_1+i f_2$, where $f_1,f_2\in V_{{\sigma}_0}\subset V {\otimes}_{\mathbb Q} \mathbb R$ are defined as follows: 
$$
f_1=\frac{1}{\sqrt{{{\sigma}_0}(d_1)}}\cdot u_1 + \sum_{i\geq 3}{\frac{x_i}{\sqrt{\mid  {{\sigma}_0}(d_i)  \mid }} \cdot u_i}, 
$$

$$
f_2=\frac{1}{\sqrt{{{\sigma}_0}(d_2)}}\cdot u_2 + \sum_{i\geq 3}{\frac{y_i}{\sqrt{\mid  {{\sigma}_0}(d_i)  \mid }} \cdot u_i}, 
$$
where $(x_3,\cdots , x_m)\in {\mathbb R}^{m-2}$ is a (very) general vector and $(y_3,\cdots , y_m)\in {\mathbb R}^{m-2}$ is a (very) general vector such that $\sum_{i=3}^{m}{x_i\cdot y_i}=0$ and $\sum_{i=3}^{m}{(x_i)^2}=\sum_{i=3}^{m}{(y_i)^2}$.\\

On multiplying $x_i$ and $y_j$ by the same (very) general $\lambda \in (0,1)$, one can guarantee that $(tr\circ \Phi)(f,f)=(tr\circ \Phi)(f_1,f_1)-(tr\circ \Phi)(f_2,f_2)=0$ and $(tr\circ \Phi)(f,\bar{f})=(tr\circ \Phi)(f_1,f_1)+(tr\circ \Phi)(f_2,f_2)=2\cdot (1- \sum_{i=3}^{m}{(x_i)^2})>0$.\\

If $u=\sum_{i}{a_i u_i}\in V$ and $(tr\circ \Phi)(u,f)=0$, then $0= (tr\circ \Phi)(u,f_1)={{\sigma}_0}(a_1)\cdot \sqrt{{{\sigma}_0}(d_1)}-\sum_{i\geq 3}{{{\sigma}_0}(a_i) \cdot \sqrt{\mid {{\sigma}_0}(d_i) \mid } \cdot x_i}$, which implies that $a_i=0$ for any $i$, since all $x_j$ are (very) general (we use the fact that $E$ is countable, while $\mathbb R$ is not). Hence $u=0$.\\

Hence $f$ gives a weight $2$ Hodge structure on $V$. By the Torelli theorem for $K3$ surfaces \cite{PS} there is a $K3$ surface $X$ and a marking $H^2(X,\mathbb Z)\cong \Lambda$, which identifies $V$ with the transcendental $\mathbb Q$-lattice of $X$, as well as their polarized Hodge structures.\\

Let us check that $E=End_{Hdg}(V)$. By a theorem of Zarhin (Theorem 1.6 in \cite{Zarhin}) $V$ is a simple $Hdg$-module and $End_{Hdg}(V)$ is a (totally real or CM-) field of algebraic numbers. Since $f\in V_{{\sigma}_0}$, we have an embedding $E\subset End_{Hdg}(V)$. Note that the restriction of the involution on $End_{Hdg}(V)$ constructed in \cite{Zarhin} to $E$ coincides with the CM-involution on $E$ (if $E$ is a CM-field) or with the identity (if $E$ is totally real).\\

Hence by Schur's lemma it is enough to show that any $\alpha \in End_{Hdg}(V)$ has an eigenvector on $V_{{\sigma}_0}\subset V{\otimes}_{\mathbb Q} \mathbb C$ with eigenvalue in ${{\sigma}_0}(E)\subset \mathbb C$. Note also that any $\alpha \in End_{Hdg}(V)$ is $E$-linear and hence ${\alpha}(V_{\sigma})\subset V_{\sigma}$. We may assume that $\alpha \neq 0$, i.e. is an $E$-linear isomorphism of $V$.\\

Let ${\alpha}(f_1+if_2)=(c_1+ic_2)\cdot (f_1+if_2)$, where $c_1, c_2 \in \mathbb R$. We need to show that $c_2=0$ and $c_1\in {{\sigma}_0}(E)\subset \mathbb R$. If we take ${\alpha}(u_i)=\sum_{j}{{{\alpha}_{ij}}\cdot u_j}$ and ${a}_{ij}={{\sigma}_0}({\alpha}_{ij})$, then 
$$
{\alpha}(f_1)=\sum_{j}{\left( \frac{a_{1j}}{\sqrt{{{\sigma}_0}(d_1)}} +\sum_{i\geq 3}{ \frac{a_{ij}}{\sqrt{\mid {{\sigma}_0}(d_i) \mid }} \cdot x_i }       \right)\cdot u_j} = c_1\cdot f_1 - c_2 \cdot f_2
$$

$$
{\alpha}(f_2)=\sum_{j}{\left( \frac{a_{2j}}{\sqrt{{{\sigma}_0}(d_2)}} +\sum_{i\geq 3}{ \frac{a_{ij}}{\sqrt{\mid {{\sigma}_0}(d_i) \mid }} \cdot y_i }       \right)\cdot u_j} = c_1\cdot f_2 + c_2 \cdot f_1
$$

This implies that $c_1 = a_{11}+ \sum_{i\geq 3}{ \frac{a_{i1}\cdot \sqrt{{{\sigma}_0}(d_1)}}{\sqrt{\mid {{\sigma}_0}(d_i) \mid }} \cdot x_i }=$  $ a_{22}+ \sum_{i\geq 3}{ \frac{a_{i2}\cdot \sqrt{{{\sigma}_0}(d_2)}}{\sqrt{\mid {{\sigma}_0}(d_i) \mid }} \cdot y_i }  $ and 

$c_2 = -a_{12}\cdot \sqrt{\frac{{{\sigma}_0}(d_2)}{{{\sigma}_0}(d_1)}} - \sum_{i\geq 3}{ \frac{a_{i2}\cdot \sqrt{{{\sigma}_0}(d_2)}}{\sqrt{\mid {{\sigma}_0}(d_i) \mid }} \cdot x_i }=$  $ a_{21}\cdot \sqrt{\frac{{{\sigma}_0}(d_1)}{{{\sigma}_0}(d_2)}} + \sum_{i\geq 3}{ \frac{a_{i1}\cdot \sqrt{{{\sigma}_0}(d_1)}}{\sqrt{\mid {{\sigma}_0}(d_i) \mid }} \cdot y_i }  $. Since $x_i$ and $y_j$ are (very) general, $a_{i1}=a_{i2}=0$ for any $i\geq 3$. Hence $c_1=a_{11}=a_{22}\in {{\sigma}_0}(E)\subset \mathbb R$ and $c_2 = -a_{12}\cdot \sqrt{\frac{{{\sigma}_0}(d_2)}{{{\sigma}_0}(d_1)}} = a_{21}\cdot \sqrt{\frac{{{\sigma}_0}(d_1)}{{{\sigma}_0}(d_2)}}$.\\

Now equalities above imply (since $x_i$ and $y_j$ are very general) that $c_2=a_{12}=a_{21}=0$. Hence $\alpha \in E$ and so $E=End_{Hdg}(V)$.\\

Case 2: $E$ is totally real and $m=3$. Let $f=f_1+i f_2$, where $f_1,f_2\in V_{{\sigma}_0}\subset V {\otimes}_{\mathbb Q} \mathbb R$ are defined as follows: 
$$
f_1=\frac{1}{\sqrt{{{\sigma}_0}(d_1)}}\cdot u_1 + \frac{x_3}{\sqrt{\mid  {{\sigma}_0}(d_3)  \mid }} \cdot u_3, 
$$

$$
f_2=\frac{\sqrt{1-(x_3)^2}}{\sqrt{{{\sigma}_0}(d_2)}}\cdot u_2, 
$$
where $x_3\in (0,1)\subset \mathbb R$ is (very) general.\\

Then $(tr\circ \Phi)(f,f)=0$, $(tr\circ \Phi)(f,\bar{f})=2(1-(x_3)^2)>0$ and if $(tr\circ \Phi)(a_1u_1+a_2u_2+a_3u_3,f)=0$ (where $a_1, a_2, a_3\in E$), then $a_1=a_2=a_3=0$.\\

If $\alpha \colon V\xrightarrow{~} V$ is a $E$-linear automorphism such that ${\alpha}(f_1+if_2)=(c_1+ic_2)\cdot (f_1+if_2)$, where $c_1, c_2 \in \mathbb R$, then one sees as in Case 1 above that $c_1\in {{\sigma}_0}(E)$, $c_2=0$ and so $\alpha$ is a left multiplication by some $a\in E$.\\

Case 3: $E$ is totally real and $m=2$. We will see that in this case a suitable $f\in V {\otimes}_{\mathbb Q} \mathbb C$ does not exist. In fact, if $f=f_1+if_2$, where $f_1, f_2 \in V_{{\sigma}_0}\subset V {\otimes}_{\mathbb Q} \mathbb R$, then we should have: 
$$
f_1=\frac{1}{\sqrt{{{\sigma}_0}(d_1)}}\cdot u_1 + \frac{x}{\sqrt{  {{\sigma}_0}(d_2)  }} \cdot u_2, 
$$

$$
f_2=-\frac{x}{\sqrt{{{\sigma}_0}(d_1)}}\cdot u_1 + \frac{1}{\sqrt{  {{\sigma}_0}(d_2)  }} \cdot u_2, 
$$
where $x\in \mathbb R$.\\

However, a $E$-linear automorphism of $V=E\cdot u_1 \oplus E\cdot u_2$ given by a matrix $\left[
   \begin{array}{ccccccccccc}
0 & -d_1 \\
d_2 & 0
   \end{array}
 \right]$ maps $f$ to $\sqrt{-{{\sigma}_0}(d_1){{\sigma}_0}(d_2)}\cdot f$ and so is an element of $End_{Hdg}(V)$, i.e. $E\neq End_{Hdg}(V)$ in this case.\\
 
Case 4: $E$ is a CM-field and $m\geq 1$. Let 
$$
f=\frac{1}{\sqrt{{{\sigma}_0}(d_1)}} \cdot u_1 {\otimes}_{E,{\sigma}_0} 1 + \sum_{i\geq 2}{ \frac{x_i}{\sqrt{ \mid {{\sigma}_0}(d_i) \mid }} \cdot u_i {\otimes}_{E,{\sigma}_0} 1 },
$$
where $(x_2, \cdots , x_m)\in {\mathbb C}^{m-1}$ is a (very) general vector such that $\sum_{i\geq 2}{\mid x_i  {\mid}^2 }\in [ 0, 1 )\subset \mathbb R$.\\

Then $(tr\circ \Phi)(f,f)=0$, $(tr\circ \Phi)(f,\bar{f})=1-\sum_{i\geq 2}{\mid x_i {\mid}^2}>0$ and if $(tr\circ \Phi)(u,f)=0$ for some $u\in V$, then $u=0$.\\

If $m=1$, then equality $E=End_{Hdg}(V)$ is automatic. So, we may assume that $m\geq 2$.\\

Let $\alpha \in End_{Hdg}(V)$ be $E$-linear, ${\alpha}(u_i)=\sum_{j}{{\alpha}_{ij} u_{ij}}$, ${\alpha}_{ij}\in E$, $a_{ij}={{\sigma}_0}({\alpha}_{ij})$ such that ${\alpha}(f)=c\cdot f$, where $c\in \mathbb C$. We have to show that $c\in {{\sigma}_0}(E)\subset \mathbb C$.\\

We have: $c=a_{11}+ \sum_{i\geq 2}{\frac{a_{i1}\cdot \sqrt{{{\sigma}_0}(d_1) } }{\sqrt{ \mid {{\sigma}_0}(d_i)  \mid  }}\cdot x_i}$ and for any $j\geq 2$, $\frac{x_j}{\sqrt{ \mid {{\sigma}_0}(d_j)  \mid  }  }\cdot c=\frac{a_{1j}}{\sqrt{{{\sigma}_0}(d_1)}}+ \sum_{i\geq 2}{\frac{a_{ij}}{\sqrt{ \mid {{\sigma}_0}(d_i)  \mid  }}\cdot x_i} $. Since $x_i$ are (very) general, this implies that $a_{i1}=0$ for all $i\geq 2$, i.e. $c=a_{11}\in {{\sigma}_0}(E)$. {\it QED}\\

{\bf Corollary 3.} {\it An $O_E$-lattice $\Phi\colon H {\otimes}_{O_E} H\rightarrow O_E$ is of $K3$ type, if and only if the following conditions hold:
\begin{enumerate}
\item[(1)] there is a primitive embedding of lattices (over $\mathbb Z$) $(H,tr\circ \Phi)\hookrightarrow \Lambda$,
\item[(2)] if $E$ is totally real, then $m\geq 3$ and there is ${\sigma}_0\in {\Sigma}_{\mathbb R}$ such that $sign ({\Phi}_{{\sigma}_0})=(2,m-2)$ and for any ${\sigma}\in {\Sigma}_{\mathbb R}$, $\sigma \neq {\sigma}_0$, $sign ({\Phi}_{{\sigma}})=(0,m)$, 
\item[(3)] if $E$ is a CM-field, then $m\geq 1$ and there is ${\sigma}_0\in {\Sigma}_{\mathbb C}$ such that $sign ({\Phi}_{{\sigma}_0})=(1,m-1)$ and for any ${\sigma}\in {\Sigma}_{\mathbb C}$, $\sigma \neq {\sigma}_0, \bar{{\sigma}_0}$, $sign ({\Phi}_{{\sigma}})=(0,m)$.
\end{enumerate}}

{\bf Remark.} (1) Necessary and sufficient conditions for the existence of a primitive embedding of two lattices were obtained by Nikulin \cite{Nikulin}. (See also \cite{Morrison} for the special case of primitive embeddings into $\Lambda$.)\\

(2) According to \cite{Witt} there exists an embedding of $\mathbb Q$-lattices $M\hookrightarrow N$, if and only if there exists a $\mathbb Q$-lattice $M'$ of signature $sign(M')=sign(N)-sign(M)$ such that $disc(M')=disc(N)\cdot disc(M)\in {{\mathbb Q}^{*}}/({\mathbb Q}^{*})^2$ and for any prime $p\in \mathbb Z$, ${\epsilon}_p(M')={\epsilon}_p(N)\cdot {\epsilon}_p(M)\cdot (disc(M),disc(M'))_p={\epsilon}_p(N)\cdot {\epsilon}_p(M)\cdot (disc(M),-disc(N))_p$ (here ${\epsilon}_p(-)$ is the Hasse invariant and $(-,-)_p$ is the Hilbert symbol). This is equivalent to the following condition:
\begin{itemize}
\item $sign(N)\geq sign(M)$ (which in particular means that $rk(M)\leq rk(N)$),
\item $(disc(M),-disc(N))_{\infty}=(-1)^{m_{-}\cdot (n_{-}+1)}$ (if we denote $sign(N)=(n_{+},n_{-})$ and $sign(M)=(m_{+},m_{-})$),
\item if $rk(M)=rk(N)-1$, then ${{\epsilon}_p}(N)\cdot {{\epsilon}_p}(M)=(disc(M),-disc(N))_p$ for all primes $p\in \mathbb Z$,
\item if $rk(M)=rk(N)-2$ and $(-disc(N)\cdot disc(M))\in ({\mathbb Q}_{p}^{*})^2$, then ${{\epsilon}_p}(N)\cdot {{\epsilon}_p}(M)=(disc(M),-1)_p$. 
\end{itemize}

In case $N={\Lambda} {\otimes}_{\mathbb Z} \mathbb Q$ we have $disc(N)=-1$, $sign(N)=(3,19)$ and ${\epsilon}_p(N)=(-1)^{p-1}$. Let us also assume that $sign(M)=(2,t)$. Then there exists an embedding of $\mathbb Q$-lattices $M\hookrightarrow {\Lambda} {\otimes}_{\mathbb Z} \mathbb Q$, if and only if:
\begin{itemize}
\item $t\leq 19$ and
\item if $t=19$ or ($t=18$ and $disc(M)\in ({\mathbb Q}_{p}^{*})^2$), then ${\epsilon}_p(M)=(-1)^{p-1}$.
\end{itemize}

(3) Diagonalizing $\Phi = d_1\cdot (X_1)^2+d_2\cdot (X_2)^2+\cdots d_m\cdot (X_m)^2$ we can decompose $\mathbb Q$-lattice $(V, tr\circ \Phi)$ as follows. Note that $\{ e_i\cdot u_j  \}$ is a $\mathbb Q$-basis of $V$ and $(tr\circ \Phi)(e_i\cdot u_j, e_{i'}\cdot u_{j'})={\delta}_{jj'}\cdot tr(e_i\cdot \overline{e_{i'}} \cdot d_j)$.\\

For any totally real number field $F$ and $a\in F$, $a\neq 0$ let us denote by $Tr_{F}<a>$ the scaled trace form on $F$, i.e. $Tr_{F}<a> \colon F {\otimes}_{\mathbb Q} F\rightarrow \mathbb Q$, $(\xi,\zeta)\mapsto Tr_F(\xi \cdot \zeta \cdot a)$. This is a nondegenerate symmetric bilinear form with discriminant $disc(Tr_{F}<a>)=Norm_F(a)\cdot d_F$, where $d_F$ is the discriminant of $F$ (see \cite{Serre}, for example).\\

If $E$ is totally real, we conclude that $(V, tr\circ \Phi)\cong \bigoplus_{j=1}^{m}Tr_F<d_j>$ as $\mathbb Q$-lattices.\\

If $E=E_0(\theta)$ is a CM-field, then for any $i,i'=1,...,s$, $(tr\circ \Phi)(e_i\cdot u_j, e_{i'+s}\cdot u_j)=tr(e_i\cdot \overline{e_{i'+s}} \cdot d_j)=tr_E(-e_i\cdot e_{i'}\cdot \theta \cdot d_j)=\sum_{\sigma \in {\Sigma}_{\mathbb C}}{(   {\sigma}(-e_i\cdot e_{i'}\cdot \theta \cdot d_j) +{\bar{\sigma}}(-e_i\cdot e_{i'}\cdot \theta \cdot d_j)  )}=0$ (we are assuming that $e_{i'+s}=\theta \cdot e_{i'}$, $\bar{\theta}=-\theta$, $\overline{e_{i'}}=e_{i'}$), $(tr\circ \Phi)(e_i\cdot u_j, e_{i'}\cdot u_j)=tr_E(e_i\cdot e_{i'}\cdot d_j)=2\cdot tr_{E_0}(e_ie_{i'}\cdot d_j)$ and $(tr\circ \Phi)(e_{i+s}\cdot u_j, e_{i'+s}\cdot u_j)=tr_E(e_i\cdot e_{i'}\cdot (-{\theta}^2 \cdot d_j))=2\cdot tr_{E_0}(e_ie_{i'}\cdot (-{\theta}^2 \cdot d_j))$. Hence in this case $(V, tr\circ \Phi)\cong \bigoplus_{j=1}^{m}(Tr_{E_0}<d_j>(2)\oplus  Tr_{E_0}<-{\theta}^2d_j>(2)  )$ as $\mathbb Q$-lattices.\\

{\bf Corollary 4.} {\it If $E$ is totally real, then $(V,\Phi)$ is of $K3$ type, if and only if the following conditions are satisfied:
\begin{enumerate}
\item[(1)] $3\leq m\leq 21$, $1\leq r\leq \frac{21}{m}$,
\item[(2)] there is ${\sigma}_0\in {\Sigma}_{\mathbb R}$ such that $sign ({\Phi}_{{\sigma}_0})=(2,m-2)$ and for any ${\sigma}\in {\Sigma}_{\mathbb R}$, $\sigma \neq {\sigma}_0$, $sign ({\Phi}_{{\sigma}})=(0,m)$,
\item[(3)] if $mr=21$ or ($mr=20$ and $(d_E)^m\cdot Norm_E(disc(\Phi))\in ({\mathbb Q}_{p}^{*})^2$), then ${\epsilon}_p(\oplus_{j=1}^{m} Tr_E<d_j> )=(-1)^{p-1}$.
\end{enumerate}}

{\it Proof:} $t=mr-2$, $disc(M)=\prod_{j=1}^{m}{disc(Tr_E<d_j>)}=(d_E)^m\cdot Norm_E(disc(\Phi))$.	{\it QED}\\

{\bf Corollary 4'.} {\it If $E=E_0(\theta)$ is a CM-field, then $(V,\Phi)$ is of $K3$ type, if and only if the following conditions are satisfied:
\begin{enumerate}
\item[(1)] $1\leq m\leq 10$, $1\leq s\leq \frac{10}{m}$,
\item[(2)] there is ${\sigma}_0\in {\Sigma}_{\mathbb C}$ such that $sign ({\Phi}_{{\sigma}_0})=(1,m-1)$ and for any ${\sigma}\in {\Sigma}_{\mathbb C}$, $\sigma \neq {\sigma}_0, \bar{{\sigma}_0}$, $sign ({\Phi}_{{\sigma}})=(0,m)$,
\item[(3)] if $ms=10$ and $(Norm_{E/{\mathbb Q}}(\theta))^m\in ({\mathbb Q}_{p}^{*})^2$, then ${\epsilon}_p(\oplus_{j=1}^{m} Tr_{E_0}<d_j> )\cdot {\epsilon}_p(\oplus_{j=1}^{m} Tr_{E_0}<-{\theta}^2\cdot d_j> )=(-1)^{p-1}\cdot ( 2^s\cdot (d_{E_0})^m\cdot Norm_{E_0}(disc(\Phi)) ,-1)_p$.
\end{enumerate}}

{\it Proof:} $t=2ms-2$, $disc(M)=\prod_{j=1}^{m}{(2^s\cdot disc(Tr_{E_0}<d_j>)\cdot 2^s\cdot disc(Tr_{E_0}<-{\theta}^2\cdot d_j>)     )}=4^{ms}\cdot (d_{E_0})^{2m}\cdot (Norm_{E_0}(-{\theta}^2))^m\cdot (Norm_{E_0}(disc(\Phi)))^2$. For any $\mathbb Q$-lattice $W$, ${\epsilon}_p(W(2))={\epsilon}_p(W)\cdot (2, (-1)^{rk(W)(rk(W)-1)/2}\cdot (disc(W))^{rk(W)-1})_p$	{\it QED}\\

{\bf Corollary 5.} {\it If $E\neq \mathbb Q$ is totally real and $(V,\Phi)$ is of $K3$ type, then 
$$
(m,r)\in \{  (3,\leq 7), (4,\leq 5), (5,\leq 4), (6,3), (6,2), (7,3), (7,2), (8,2), (9,2), (10,2)  \}.
$$

If $E=\mathbb Q$ and $(V,\Phi)$ is of $K3$ type, then $3\leq m \leq 21$.\\

If $E\neq \mathbb Q$ is totally real, ( there is ${\sigma}_0\in {\Sigma}_{\mathbb R}$ such that $sign ({\Phi}_{{\sigma}_0})=(2,m-2)$ and for any ${\sigma}\in {\Sigma}_{\mathbb R}$, $\sigma \neq {\sigma}_0$, $sign ({\Phi}_{{\sigma}})=(0,m)$ ) and 
$$
(m,r)\in \{  (3,\leq 6), (4,\leq 4), (5,3), (5,2), (6,3), (6,2), (7,2), (8,2), (9,2) \},
$$
then $(V,\Phi)$ is of $K3$ type.\\

If $E= \mathbb Q$, $sign ({\Phi})=(2,m-2)$ and $3\leq m=dim(V)\leq 19$, then $(V,\Phi)$ is of $K3$ type.}\\

{\bf Corollary 5'.} {\it If $E=E_0(\theta)$ is a CM-field, ${\theta}^2\in E_0$, $E_0\neq \mathbb Q$ is the totally real subfield (fixed by the CM-involution) and $(V,\Phi)$ is of $K3$ type, then 
$$
(m,r)\in \{  (1,\leq 10), (2,\leq 5), (3,3), (3,2), (4,2), (5,2) \}.
$$

If $E=\mathbb Q (\theta)$ is a CM-field (${\theta}^2\in\mathbb Q$) and $(V,\Phi)$ is of $K3$ type, then $1\leq m \leq 10$.\\

If $E=E_0(\theta)$ is a CM-field, ${\theta}^2\in E_0$, $E_0\neq \mathbb Q$ is the totally real subfield (fixed by the CM-involution), ( there is ${\sigma}_0\in {\Sigma}_{\mathbb C}$ such that $sign ({\Phi}_{{\sigma}_0})=(1,m-1)$ and for any ${\sigma}\in {\Sigma}_{\mathbb C}$, $\sigma \neq {\sigma}_0, \bar{{\sigma}_0}$, $sign ({\Phi}_{{\sigma}})=(0,m)$ ) and 
$$
(m,r)\in \{  (1,\leq 9), (2,\leq 4), (3,3), (3,2), (4,2) \},
$$
then $(V,\Phi)$ is of $K3$ type.\\

If $E= \mathbb Q (\theta)$ is a CM-field (${\theta}^2\in\mathbb Q$), $sign ({\Phi})=(1,m-1)$ and $1\leq m=dim(V)\leq 9$, then $(V,\Phi)$ is of $K3$ type.}\\

\subsection*{Acknowledgement}

We thank Yuri Zarhin for suggesting this problem. Possible fields $E$ which can appear in the totally real case were also studied independently by Bert van Geemen \cite{vanGeemen2}.\\

\bibliographystyle{ams-plain}

\bibliography{HermitianFormsOfK3type}

\end{document}